%% file: main.tex
\documentclass[12pt]{article}
\usepackage{amsmath,amsfonts,amsthm,amssymb}
\usepackage{graphicx,epic,eepic,times}

\usepackage[margin=1in]{geometry}
\usepackage[bottom,none]{draftcopy}
\draftcopySetGrey{0.5}

\newcounter{kevinslistcountertoo}%
\renewenvironment{enumerate}{
\begin{list}{(\roman{kevinslistcountertoo})}
{\usecounter{kevinslistcountertoo}
\setlength{\parsep}{3pt}
\setlength{\labelwidth}{24pt}
\setlength{\itemsep}{3pt}
\setlength{\topsep}{3pt}}}{\end{list}}

\usepackage[small]{caption}

\newtheorem*{algorithm*}{Algorithm}

\newcommand{\cov}{{\rm Cov}}
\newcommand{\var}{{\rm Var}}
\newcommand{\E}{{\mathbb E}}

\newcommand{\prob}{{\rm Prob}}

\newcommand{\R}{{\mathbb R}}

\newcommand{\eps}{\varepsilon}

\begin{document}

\title{\Large Coupling Control Variates for Markov Chain Monte Carlo}

\author{\normalsize Jonathan B. Goodman\thanks{Courant Institute of
    Mathematical Sciences, New York University, 251 Mercer Street, New
    York, NY 10012, USA.  {\bf goodman@cims.nyu.edu}}
  {\ and} Kevin K. Lin\thanks{Department of
    Mathematics, University of Arizona, 617 N. Santa Rita Ave., Tucson,
    AZ 85721, USA.  {\bf klin@math.arizona.edu}}}

\date{\small June 6, 2008}

\maketitle

\begin{abstract}
We show that Markov couplings can be used to improve the accuracy of
Markov chain Monte Carlo calculations in some situations where the
steady-state probability distribution is not explicitly known.  The
technique generalizes the notion of control variates from classical
Monte Carlo integration.  We illustrate it using two models of
nonequilibrium transport.
\end{abstract}

%%%%%%%%%%%%%%%%%%%%%%%%%%%%%%%%%%%%%%%%%%%
\section*{Introduction}

Markov chain Monte Carlo (MCMC) algorithms generate samples from a
probability distribution by simulating a Markov chain that leaves the
distribution invariant.  One estimates expected values by time averaging
over long simulations~\cite{KW,sokal}.  For high-accuracy Monte Carlo
computations, variance reduction methods are crucial.  Unfortunately,
some variance reduction methods are hard to apply in MCMC, particularly
when there is no explicit expression for the steady-state probability
distribution of the Markov chain.

In this paper, we demonstrate a technique for MCMC variance reduction
which can improve accuracy by factors of up to 2 or more in certain
situations where an {\em approximate} steady-state distribution is
known.  The technique, which we call {\em coupling control variates},
builds on earlier work using Markov couplings in MCMC~\cite{neal,o1,o2}.
Specifically, we assume that we can obtain an explicit approximation of
the steady-state distribution, and that the expected values of this
approximate distribution are known.  The basic idea is to find a second
Markov process which (i) leaves the approximate distribution invariant,
and (ii) ``shadows'' ({\em i.e.}, closely follows) the original Markov
process.  The expectations of the approximate distribution then provide an
initial ``guess,'' which we correct by simulating the two ``coupled''
processes to estimate the difference (in expected values) between the
true steady-state distribution and our approximate distribution.

We apply the technique to certain lattice models from statistical
physics, in which the steady-state probability distribution is
approximately a product of local distributions when the system is out of
equilibrium.\footnote{Here, ``equilibrium'' is used in the sense of
  statistical physics, {\em i.e.}, ``thermal equilibrium.'' This means
  that the Markov chain satisfies detailed balance~\cite{KW}, and the
  steady-state probability distribution is a Gibbs-Boltzmann
  distribution $\frac{1}{Z}e^{-\beta H}$.  Steady-state distributions of
  Markov chains that are {\em not} in equilibrium are known as
  ``nonequilibrium steady states.''  We focus on the latter here.}
These systems are of interest in the theory of transport processes such
as heat conduction.  In this paper, we consider models consisting of a
linear chain of lattice sites coupled to ``heat baths'' at each end;
each bath is characterized by thermodynamic parameter(s) like
temperature, chemical potential, {\em etc.}  The steady-state
probability distribution is a Gibbs-Boltzmann distribution if the bath
parameters are equal.  This is not the case for unequal heat baths.
However, a large lattice out of equilibrium may still have a
steady-state distribution that is {\em locally} in equilibrium, {\em
  e.g.}, for heat flow, the statistics at a given location is
approximately governed by a Gibbs-Boltzmann distribution with a local
temperature (see Sect.~\ref{sect2: applications} for details).  We will
show how such ``local equilibrium'' distributions can be used to achieve
variance reduction.

We note that similar ideas have appeared in previous studies.  In
addition to the works cited earlier, another example is ``shadow hybrid
Monte Carlo'' in molecular dynamics~\cite{IH}.  See also {\cite{glynn}}
for a version of this idea applied to Markov sensitivity analysis.
Finally, we point out that Markov couplings have been used in a quite
different way to perform exact Monte Carlo sampling~\cite{propp}.

%%%%%%%%%%%%%%%%%%%%%%%%%%%%%%%%%%%%%%%%%%%
\section{Coupling control variates}
\label{sect1: coupling control variates}

\subsection{General framework}
\label{General framework}

We begin by recalling the technique of control variates in classical
Monte Carlo (MC) integration~\cite{HH}: suppose $X$ is a random variable
with probability density $p_X$, and we want to estimate its expected
value $\bar{X} = \E[X] = \int_{-\infty}^\infty{x\cdot p_X(x)~dx}$.  The
standard Monte Carlo estimator of $\bar{X}$ is
\begin{equation}
\widehat{X}_n = \frac{1}{n}\sum_{k=1}^n X_k~,
\label{se}
\end{equation}
where $X_1, X_2, \cdots,$ are independent samples from the distribution
$p_X$.  The variance of the estimator is $\var[\widehat{X}_n] =
\var[X]/n$.  It is not generally possible to improve the $c/n$ scaling;
more accurate estimates are usually obtained by reducing the variance of
the estimand.

A {\em control variate} for $X$ is a random variable $Y$ whose expected
value $\bar Y=\E[Y]$ is known and is correlated with $X$.  One can
estimate $\bar{X}$ using the {\em control variate estimator}
\begin{equation}
\widehat{X}_{CV,\alpha;n} = \frac{1}{n}\sum_{k=1}^n \big[ X_k +
  \alpha\cdot(\bar{Y} - Y_k) \big]\ .
\label{cve}
\end{equation}
where $(X_k,Y_k), k=1, 2, \cdots,$ are samples from the joint
distribution of $X$ and $Y$, and $\alpha$ is an adjustable parameter.
Optimizing $\var[\widehat{X}_{CV,\alpha;n}]$ over $\alpha$ gives an
optimal control variate estimator of $\bar{X}$ with variance
$$
\frac{1}{n}~\var[X] \cdot ( 1 - \rho^2_{XY} )~,
$$ where $\rho_{XY}$ is the correlation coefficient
$\cov(X,Y)/(\var[X]\cdot\var[Y])^{1/2}$.  In the special case
$\alpha=1$, Eq.~(\ref{cve}) simply corrects the initial ``guess''
$\bar{Y}$ with an estimate of $\bar{X}-\bar{Y}$.

%% Control variates are particularly effective in cases that have
%% approximations with known expected values.  In the present case,
%% the steady-state of a long chain should have an approximate product
%% form.

\bigskip
Consider now {\em Markov chain Monte Carlo}, where the samples are not
independent, but are successive states of a Markov process.  For
concreteness, let $X_t$ be a time-homogeneous continuous-time Markov
process with finite state space\footnote{Extending our ideas to more
  general settings is straightforward.  See for instance Sect.~\ref{Sect
    KMP model}.}  $\Omega$.  The dynamics of $X_t$ are completely
specified by the {\em transition rates} $R(x'|x)$, which tell us the
rate at which $X_t$ jumps from state $x$ to state $x'$, {\em i.e.},
$\prob\Big(X_{t+\Delta t} = x' \Big| X_t = x\Big) = R(x'|x)\cdot\Delta t
+ O(\Delta t^2)~.$ We assume that the process $X_t$ has a unique
steady-state probability distribution $P$, so that
$\sum_{x'}R(x|x')~P(x') = \sum_{x'}R(x'|x)~P(x)$.

Given an observable $\phi:\Omega\to\R$, one can obtain a direct estimate
of $\E_X[\phi] = \sum_{x\in\Omega} \phi(x)\cdot P(x)$ by simulating the
process $X_t$ for $t\in[0,T]$ and applying the simple estimator
\begin{equation}
\widehat{\phi}_T = \frac{1}{T}\int_0^T{\phi(X_t)~dt}~.
\label{simple estimator}
\end{equation}
This converges almost surely to $\E_X[\phi]$ as $T\to\infty$.  The
variance of $\widehat{\phi}_T$ is given by the {\em Kubo variance
  formula}~\cite{anderson}
\begin{equation}
\frac{\var[\phi]\cdot\tau}{T} + O(1/T^2)~,
\label{kubo}
\end{equation}
where $\var[\phi]$ is the variance of the observable $\phi$ with respect
to $P$.  The constant $\tau$ is the {\em integrated autocorrelation
  time}
$$
\tau = \int_{-\infty}^{\infty} \rho(t)\ dt~,
$$ where $\rho(t) =C(t)/C(0)$ is the {\em time-autocorrelation function}
of $\phi(X_t)$, and
$$
C(t) = \lim_{t_0\to\infty}\cov\Big(\phi(X_{t+t_0}), \phi(X_{t_0})\Big)~.
$$ Note that $\tau$ depends on both the observable $\phi$ and the Markov
process $X_t$.

As in the case of MC integration, it is not generally possible to
improve the $c/T$ scaling in Eq.~(\ref{kubo}).  Variance reduction
schemes typically aim to reduce either the autocorrelation time $\tau$
or the variance $C(0)$ of the estimand.

To extend the notion of control variates to this setting, one looks for
a second Markov process $Y_t$ which is correlated to the process of
interest $X_t$~\cite{neal,o1,o2}.  The notion of correlated processes
can be made precise by {\it Markov couplings}~\cite{lindvall}: if $X_t$
and $Y_t$ are Markov processes with respective transition rates $R_X$
and $R_Y$, a Markov coupling of $X_t$ and $Y_t$ is a specification of
{\em joint transition rates} $R_{XY}((x',y')|(x,y))$ for transitions
from $(X_t,Y_t)=(x,y)$ to $(X_t,Y_t)=(x',y')$, so that
\begin{equation}
  \begin{array}{l}
    \sum_{y'}R_{XY}((x',y')|(x,y)) = R_X(x'|x)\mbox{~~~for all
      $y, x, x'$~,\qquad and}\\[1ex]
    \sum_{x'}R_{XY}((x',y')|(x,y)) = R_Y(y'|y)\mbox{~~~for all $x, y, y'$}~.
  \end{array}
  \label{coupling condition}
\end{equation}
In other words, a Markov coupling of $X_t$ and $Y_t$ is a Markov process
on the product space $\Omega\times\Omega$ that gives a realization of
$X_t$ when projected onto the first component, and likewise gives $Y_t$
when projected onto the second.

Suppose a process $Y_t$ can be found such that the expectation
$\E_Y[\phi]$ with respect to the stationary distribution $Q$ of $Y_t$
can be computed easily.  We define the {\em coupling control variate
  estimator} by
\begin{equation}
\widehat{\phi}_{couple,\alpha} = \frac{1}{T}\int_0^T \Big[ \phi(X_t) +
  \alpha\cdot\big(\E_Y[\phi] - \phi(Y_t)\big)\Big]~dt~.
\label{ccv}
\end{equation}
The process $Y_t$ is the {\em coupling control variate}.  It is possible
to estimate a nearly optimal $\alpha$ using the Kubo variance formula
(\ref{kubo}), but for simplicity we will always set $\alpha=1$ in this
paper.\footnote{For the models studied in this paper, it is expected
  that the optimal $\alpha$ will be $\approx 1$.  In more general
  situations, it is important (and not difficult) to estimate an optimal
  $\alpha$.}  In order for the coupling control variate to be effective
with this choice of $\alpha$, $\phi(Y_t)-\phi(X_t)$ should have small
variance, {\em i.e.}, the states $X_t$ and $Y_t$ should remain as close
to each other as possible.

\subsection{The coupling control variate algorithm}
\label{metropolization}

Now, suppose we are interested in computing $\E_X[\phi]$ for a Markov
process $X_t$ with transition rates $R_X(x'|x)$.  Suppose further that
the steady-state distribution $P$ is not known, but that an approximate
steady-state distribution $Q$ is available.  Our aim is to construct a
coupled process $(X_t,Y_t)$ with transition rates
$R_{XY}((x',y')|(x,y))$ so that
\begin{enumerate}

  \item The marginal $X_t$ has transition rates $R_X$, and therefore
    steady-state distribution $P$.

  \item The marginal $Y_t$ has steady-state distribution $Q$.

  \item $X_t$ and $Y_t$ remain as close as possible given constraints
    (i) and (ii).

\end{enumerate}
We show here how the coupling $R_{XY}((x',y')|(x,y))$ can be constructed
from a coupling $R_{XX}$ of two realizations of $R_X$ processes.  Such
couplings are available in many situations; see Sect.~\ref{some
  practical considerations}.  The basic idea is to apply the
Metropolis-Hastings algorithm using the second component of $R_{XX}$ as
proposal and the distribution $Q$ as the target distribution.  The
result is a process $Y_t$ satisfying the detailed balance condition
with respect to $Q$:
\begin{equation}
  Q(y')\cdot R_Y(y|y') = Q(y)\cdot R_Y(y'|y)~.
  \label{detailed balance}
\end{equation}
Thus, the stationary distribution of $Y_t$ is $Q$.  This is a
straightforward generalization of the detailed balance condition for
discrete time Markov chains; see, {\em e.g.},~\cite{KW,sokal}.

More precisely, recall that one way to simulate continuous-time
finite-state Markov processes is as follows (sometimes known as the
Gillespie algorithm~\cite{gillespie}): let $R(x) = \sum_{x'\neq
  x}R(x'|x)$ be the {\em total exit rate} from a state $x\in\Omega$.
Let $T_n$ be the times at which the system jumps to the next state, and
let $X(n) = X_{T_n+}$ be the state of the system after each jump.  If
$X(n)=x$, we set an exponential clock of mean $1/R(x)$.  When the clock
rings, we choose a new state $x'$ with probability $P(x'|x) =
R(x'|x)/R(x)$ and set $X(n+1) = x'$.  Note that $X_t = X(n)$ for
$T_n\leq t < T_{n+1}$.

The following simple algorithm generates one step of a coupled process
$(X_t,Y_t)$ satisfying conditions (i-iii) above:

\medskip
\begin{algorithm*}
Let ${\bf State} = (x,y)$ be the current state of the joint process
$(X_t,Y_t)$.  With rate $R_{XX}(x',y'|x,y)$, set ${\bf Proposal} =
(x',y')$.  Compute
  \begin{equation}
    Z =\frac{Q(x')\cdot R_X(x|x')}{Q(x)\cdot R_X(x'|x)}~.
    \label{eq metropolis ratio}
  \end{equation}
  \begin{itemize}

  \item[] With probability $\min(Z,1)$, we accept ${\bf Proposal}$ and
    set {\bf NewState} to $(x',y')$.

  \item[] With probability $1-\min(Z,1)$, we reject ${\bf Proposal}$ and
    set {\bf NewState} to $(x',y)$.

  \end{itemize}
\end{algorithm*}

\medskip
\noindent
It is easy to check that the coupled process $(X_t,Y_t)$ generated by
this algorithm satisfies Eq.~(\ref{detailed balance}).  Thus, the
estimator (\ref{ccv}), when applied to $(X_t,Y_t)$, is always consistent
in that $\widehat{\phi}_{couple;T}\to\E_Y[\phi]$ as $T\to\infty$.  Note,
however, that whether the variance of the coupling control variate
estimator is lower than that of the simple estimator (\ref{simple
  estimator}) depends on the coupling $R_{XX}$ and the approximate
distribution $Q$.

\paragraph{Remark.}
We note that when computing the expectation of static observables using
this algorithm for continuous-time Markov chains, one can reduce
variance a little bit more by replacing the time intervals $T_{n+1}-T_n$
by the mean $1/R(X(n))$.

\subsection{Some practical considerations}
\label{some practical considerations}

\medskip
\noindent
{\em Approximate stationary distribution.}  The choice of $Q$ is
problem-dependent.  In the nonequilibrium models discussed in
Sect.~\ref{sect2: applications}, as in many other physical situations,
perturbative analysis of the relevant master equation often gives good
candidates for $Q$.  Note that because the coupling estimator is always
consistent, it is not necessary to know {\em a priori} how good an
approximation $Q$ is to the true stationary distribution, so that one
can take advantage of uncontrolled approximations.  However, the degree
of variance reduction depends on the distribution $Q$ and the coupling
$R_{XX}$.

To choose the distribution $Q$, one should follow these criteria:
\begin{enumerate}

\item The expected value $\E_Q[\phi]$ should be easy to compute.  This
  is necessary in order to apply the coupling control variate estimator
  (\ref{ccv}).

\item The distribution $Q$ should be ``close enough'' to the true
  stationary distribution $P_X$ that the rejection rate is low.  We may
  then expect $Y_t$ to remain close to $X_t$, so that the coupling
  control variate estimator may have low variance.

\end{enumerate}

\medskip
\noindent
{\em Constructing couplings.}  How do we obtain a coupling $R_{XX}$ to
start with?  As mentioned earlier, constructing Markov couplings is not
always straightforward.  However, couplings have long been used as a
theoretical tool for studying the ergodic properties of Markov
processes, and ``good'' couplings have been found for a broad range of
stochastic models~\cite{lindvall}.  In many (though not all) cases, it
suffices to simply use the same sequence of random numbers to couple two
Markov processes.  Examples include stochastic differential equations
that are contractive in the sense that their largest Lyapunov exponent
is negative~\cite{lejan} and the models in Sect.~\ref{sect2:
  applications}.

\medskip
\noindent
{\em Factors affecting scaling of errors.}  The variance of the coupling
control variate estimate is
\begin{equation}
\var\Big(\widehat{A}_{couple}\Big)
= \frac{\var[\phi(X)-\phi(Y)]\cdot\tau_{couple}}{T} + O(1/T^2)~,
\label{ccv error}
\end{equation}
where $\tau_{couple}$ is here the integrated autocorrelation time of
$\phi(X_t)-\phi(Y_t)$, and $\var[\phi(X)-\phi(Y)]$ is the variance of
the random variable $\phi(X)-\phi(Y)$ with respect to the stationary
distribution of the coupled process on the product space
$\Omega\times\Omega$.  Note that if the coupling is effective in keeping
$\phi(X_t)-\phi(Y_t)$ small, then the variance in Eq.~(\ref{ccv error})
will be small.  However, when a proposed move is rejected by our
algorithm, the process $Y_t$ ``stands still.''  The process $Y_t$ (and
hence $\phi(X_t)-\phi(Y_t)$) may therefore have a slower correlation
time than $X_t$.  That is, the amount by which the variance of the
estimator is reduced may reflect competition between lower variance and
larger correlation time.

\medskip
\noindent
{\em Overhead and running time.}  Another practical consideration is the
complexity of $Q$ and the coupling $R_{XX}$: a ``good'' coupling that is
computationally expensive to implement may not, in the end, be worth the
effort.  Couplings that are easy to implement, for example simply using
the same sequence of random numbers, have a distinct advantage in this
regard.

%%%%%%%%%%%%%%%%%%%%%%%%%%%%%%%%%%%%%%%%%%%
\section{Nonequilibrium transport processes}
\label{sect2: applications}

\subsection{Symmetric simple exclusion process}

%%%%%%%%%%%%%%%%%%%%%%%%%%%%%%%%%%%%%%%%%%%%%%%%%%%%%%%%%%%%%%%%%%%
\begin{figure}
\begin{center}
\includegraphics[scale=0.7,bb=0 0 220 66]{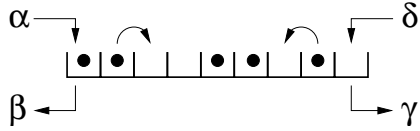}
\end{center}
\caption{The symmetric simple exclusion process.}
\label{fig ssep}
\end{figure}
%%%%%%%%%%%%%%%%%%%%%%%%%%%%%%%%%%%%%%%%%%%%%%%%%%%%%%%%%%%%%%%%%%%

The first model we consider is the {\em symmetric simple exclusion
  process} (SSEP) in one space dimension~\cite{liggett}.  This is a
stochastic lattice gas model of a linear medium with a reservoir placed
at each end.  The two reservoirs are typically maintained at different
densities, so that there is a net flow of particles through the medium.
More precisely, the domain is a linear chain of $N$ sites, with each
site holding at most one particle at any given time.  Thus, the state of
the system $\sigma\in\Omega$ can be thought of as a binary string of
length $N$, with $|\Omega| = 2^N$.  The dynamics are as follows: each
particle carries an exponential clock of rate 1.  When the clock rings,
the particle will try to jump to a neighboring site, choosing left and
right with equal probability; the particle does not move if the target
site is occupied.  The left reservoir will place a particle in site 1,
when it is unoccupied, at rate $\alpha$; and remove a particle from site
1, when it is occupied, at rate $\beta$.  The right reservoir acts on
site $N$ in an analogous manner, at rates $\delta$ and $\gamma$,
respectively.  See Fig.~\ref{fig ssep}.  Note that the total particle
number is conserved, except when the reservoirs inject or remove a
particle.

We begin by summarizing some known results on the SSEP;
see~\cite{derrida,liggett} for details.  It is easy to show that the
SSEP has a unique stationary distribution $P_N$.  Much is known about
$P_N$.  In particular, various probabilities can be calculated exactly
using the ``matrix method.''  The SSEP thus provides a convenient test
case for illustrating coupling control variates in nonequilibrium
transport models.  A central motivation for studying models like the
SSEP is to understand how macroscopic transport processes arise from
microscopic dynamics.  One quantity of interest is the macroscopic
density profile $\rho:(0,1)\to\R$, defined by
\begin{equation}
\rho(x) = \lim_{N\to\infty}\E_N[\sigma_{[xN]}] ,\qquad x\in(0,1)~,
\end{equation}
where $\E_N[\cdot]$ denotes expectation with respect to $P_N$.  Another
quantity of great interest is the correlation between distant sites (see
below).

Specifically, let $\rho_L = \alpha/(\alpha+\beta)$ and $\rho_R =
\delta/(\delta+\gamma)$.  These quantities can be thought of as the
particle densities of the reservoirs.  When $\rho_L = \rho_R = \rho_0$,
the SSEP satisfies detailed balance, and it is easy to check that the
equilibrium distribution is
\begin{equation}
P_N(\sigma) = \prod_{i=1}^N p(\sigma_i)~,
\end{equation}
where $p(1)=\rho_0$ and $p(0)=1-\rho_0$.  The occupation numbers become
IID Bernoulli random variables.  Note that this means
$\rho(x)\equiv\rho_0$.

If $\rho_L \neq \rho_R$, it can be shown that
\begin{equation}
\rho(x) = \rho_L\cdot(1-x) + \rho_R\cdot x\ .
\label{ssep profile}
\end{equation}
The non-constant profile reflects the presence of a nonzero current.
The stationary distribution $P_N$ is no longer a product: the covariance
$\cov_N(\sigma_i,\sigma_j)$ is nonzero for $i\neq j$.  The dynamics no
longer satisfies detailed balance.

The large-$N$ scaling of spatial correlations is also known.  Fix $x$,
$y$, so that $0<x<y<1$.  Then~\cite{derrida}
\begin{equation}
\lim_{N\to\infty}N\cdot\cov_N\big(\sigma_{[xN]},\sigma_{[yN]}\big)
 = -(\rho_R-\rho_L)^2\cdot x~(1-y)\ .
\end{equation}
Thus, for $N\gg 1$ and $i$, $j$ not too near the end points of $(0,1)$,
we have $\cov_N\big(\sigma_i,\sigma_j) = O(1/N)$.  We note that this
$1/N$ scaling is not unique to the SSEP --- it has been observed in
other settings as well~\cite{bertini,derrida,LY,spohn}.  The correlation
is thus quite weak for large $N$.  This means that computing
correlations in nonequilibrium transport models like the SSEP presents
numerical difficulties: when the covariances are $O(1/N)$ and the
occupation numbers $\sigma_i$ themselves remain $O(1)$, a direct
computation entails subtracting two quantities of like magnitude to
estimate a much smaller number.

\bigskip
To apply coupling control variates to this problem, we need an
approximate stationary distribution $Q$ and a coupling.  For
nonequilibrium transport models like the SSEP, a choice of $Q$ is
suggested by the notion of {\em local thermal equilibrium} (LTE): in
physical terms, even though the system cannot be in thermal equilibrium
because the two ends are in contact with reservoirs at different
densities, for large $N$ it is generally expected that small parts of
the medium will reach approximate local thermal
equilibrium~\cite{degroot}.  For the SSEP, it has been shown that LTE
holds in the following sense: fix $x\in(0,1)$ and a positive integer
$k$.  Then, as $N\to\infty$ with $x$ and $k$ fixed, the occupation
numbers $\sigma_{[xN]}, \sigma_{[xN]+1}, \cdots, \sigma_{[xN]+k}$
converge in distribution to independent, identically-distributed
Bernoulli random variables with $\prob(\sigma=1) = \rho(x)$, where
$\rho$ is the linear profile given in Eq.~(\ref{ssep profile}).
Heuristically, this tells us that even though the system cannot attain a
global thermal equilibrium when $\rho_L\neq\rho_R$, it does approach
local equilibrium when $N\gg 1$.  It also suggests that we use as our
approximate stationary distribution
\begin{equation}
Q_N(\sigma) = \prod_{i=1}^Nq_i(\sigma_i)\ ,
\end{equation}
where $q_i(1) = \rho(x_i)$, $q_i(0) = 1-\rho(x_i)$, and $x_i
=\frac{i}{N+1}$.  The distribution $Q_N$ can be thought of as a local
equilibrium distribution, in which the sites are occupied independently
with probability $\rho(x_i)$.  The LTE property suggests that $Q_N$ may
become a better approximation of $P_N$ as $N\to\infty$, at least
locally.

The other ingredient we need is $R_{XX}$, a coupling of the SSEP to
itself, so that we can use the algorithm in Sect.~\ref{metropolization}
to construct a coupling control variate.  This is
straightforward~\cite{liggett}: given two copies of SSEP, we simply
carry out the same moves in both copies whenever possible, and move
independently when not.  More precisely, let $Moves(\sigma)$ denote the
set of all available moves for $\sigma$, where a move means a particle
jumping from site $i$ to site $j$ (for all $i$, $j$ with $|i-j|=1$) or
changing the occupation number of site 1 or site $N$.  To each move in
$Moves(\sigma)\cup Moves(\tilde\sigma)$, we attach an independent
exponential clock of the appropriate rate --- 1/2 for jumps, $\alpha$
for injection by the left reservoir, {\em etc.}  When a clock goes off,
check if the corresponding move is in $Moves(\sigma)\cap
Moves(\tilde\sigma)$, {\em i.e.}, whether $\sigma$ and $\tilde\sigma$
can make the same move.  If so, update both $\sigma$ and $\tilde\sigma$
accordingly.  If the move is in $Moves(\sigma)\setminus
Moves(\tilde\sigma)$, {\em i.e.}, if only $\sigma$ can make the move,
then update only $\sigma$.  Similarly for moves in
$Moves(\tilde\sigma)\setminus Moves(\sigma)$.  This algorithm couples
two copies of the SSEP process.

%% the joint state of the coupled process, and again split the generator
%% into $L_{\mbox{\tiny SSEP}} = L_{\mbox{\tiny int}} + L_1 + L_n$:
%% \begin{align*}
%% L_{\mbox{\tiny int}}f(\sigma,\tilde\sigma)
%% &=\sum_{(i,j)\in A(\sigma,\tilde\sigma)}\big[
%% f(\sigma^{ij},\tilde\sigma^{ij})-f(\sigma,\tilde\sigma)\big]
%% +\sum_{(i,j)\notin A(\sigma,\tilde\sigma)}\big[
%% f(\sigma^{ij},\tilde\sigma) + f(\sigma,\tilde\sigma^{ij}) -
%% 2f(\sigma,\tilde\sigma)\big]\ ,\\
%% A(\sigma,\tilde\sigma) &= \{(i,j) : \sigma_i=\tilde\sigma_i,
%% \sigma_j=\tilde\sigma_j, |i-j|=1\}\ ,
%% \end{align*}
%% and the left boundary term is given by
%% \begin{equation}
%% L_1f(\sigma,\tilde\sigma)
%% =\left\{\begin{array}{ll}
%% \alpha f(\sigma^{1+},\tilde\sigma^{1+}) +
%%   \beta f(\sigma^{1-},\tilde\sigma^{1-}) -2f(\sigma,\tilde\sigma), &
%%   \sigma_1=\tilde\sigma_1=0\\
%% \alpha[f(\sigma^{1+},\tilde\sigma) + f(\sigma,\tilde\sigma^{1+}) ] +&\\
%% \beta[f(\sigma^{1-},\tilde\sigma) + f(\sigma,\tilde\sigma^{1-})] -
%% 4f(\sigma,\tilde\sigma), &\mbox{ otherwise}\\
%% \end{array}\right.
%% \end{equation}
%% The right boundary term $L_n$ is similarly defined.

We can now apply the Metropolis-Hastings construction from
Sect.~\ref{metropolization}.  This yields a coupling control variate for
the SSEP, with Metropolis ratios $Z$ given by the following table:
\begin{center}
\begin{tabular}{rll}
Transition from site $i$ to $j$, $|i-j|=1$
& $Z_{ij} = \frac{1 - \rho_i}{\rho_i}\cdot\frac{\rho_j}{1 - \rho_j}$&\\
Injection (removal) by left reservoir
& $Z_{\mbox{\tiny L,in}}=\frac{\rho_1}{1-\rho_1}\cdot\frac{\beta}{\alpha}$
& $\big(Z_{\mbox{\tiny L,out}} = 1/Z_{\mbox{\tiny L,in}}\big)$\\
Injection (removal) by right reservoir
& $Z_{\mbox{\tiny R,in}} =\frac{\rho_n}{1-\rho_n}\cdot\frac{\gamma}{\delta}$
& $\big(Z_{\mbox{\tiny R,out}} = 1/Z_{\mbox{\tiny R,in}}\big)$\\
\end{tabular}
\end{center}
Note that the $Z$ ratios involve only local quantities because the
distribution $Q_N$ has product form.  Note also that the rejection
probabilities are quite small when $N\gg 1$: since $\rho_i-\rho_j =
O(1/N)$, the Metropolis-Hastings ratios $Z$ above are $1+O(1/N)$ (as
long as $0 < \rho_L, \rho_R < 1$).  Thus, the Metropolis-Hastings
algorithm rejects fewer and fewer samples as $N\to\infty$.

\bigskip
\noindent
{\em Numerical results.}  To assess the effectiveness of the coupling
control variate, we use a metric we call the error ratio
\begin{equation}
  e_N[\phi] =
  \left(\frac{\var_N\big[\widehat{\phi}_{couple}\big]}
       {\var_N\big[\widehat{\phi}\big]}\right)^{1/2}
  \label{error ratio}
\end{equation}
for a given observable $\phi$.  The error ratio measures the amount by
which the estimator $\widehat{\phi}_{couple}$ improves the accuracy of the
estimate.

\begin{figure}
\begin{center}
\begin{tabular}{cc}
\includegraphics[bb=0in 0in 3in 2.25in]{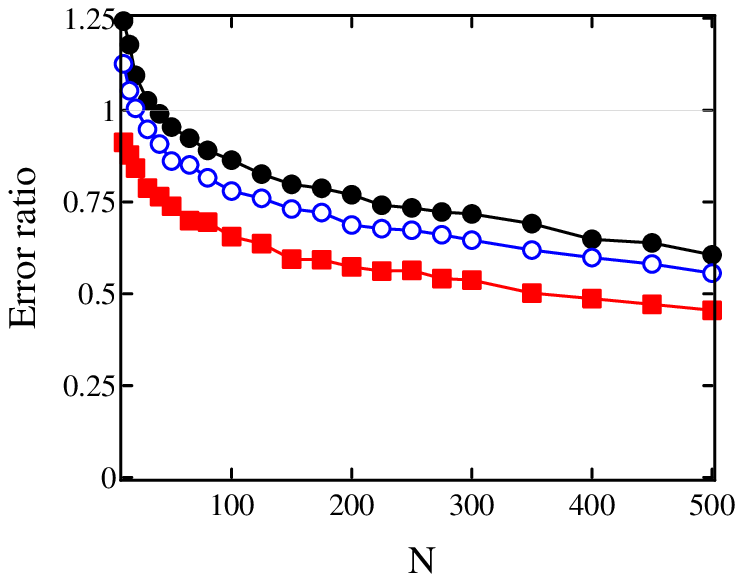}&
\includegraphics[bb=0in 0in 3in 2.25in]{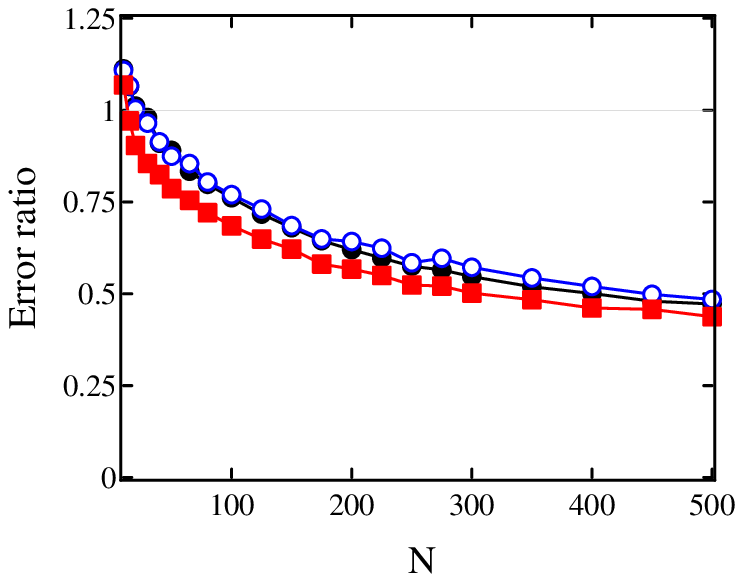}\\
(a) & (b)\\
\end{tabular}
\end{center}
\caption{The SSEP error ratio {\em vs.} system size $N$.  In (a), we
  show the error ratio $e_N$ (see text) for the estimated density at
  $x=0.3$ (solid discs), $x=0.5$ (open circles), and $x=0.8$ (squares).
  In (b), we show the error ratio for the near-neighbor product
  $\sigma_{[xN]}\cdot\sigma_{[xN]+1}$ with $x=0.5$ (solid discs), and
  for the products $\sigma_{[xN]}\cdot\sigma_{[yN]}$ with
  $(x,y)=(0.4,0.7)$ (open circles) and $(x,y)=(0.2,0.9)$ (squares).  The
  errors are estimated using batched means estimators~\cite{sokal}.  The
  parameters are $\alpha=2$, $\beta=0.1$, $\delta=0.3$, and $\gamma=1$.}
\label{fig ssep error ratio}
\end{figure}

Fig.~\ref{fig ssep error ratio}(a) shows the error ratio
$e[\sigma_{[xN]}]$ for the occupation numbers at a few selected
locations along the chain, specifically $x\in\{0.3, 0.5, 0.8\}$.  The
error ratio decreases with increasing $N$.  The improvement with $N$ is
expected, since the local equilibrium distribution $Q_N$ is expected to
be a better approximation of the true stationary distribution $P_N$ when
$N$ is big.  Indeed, our data show that the rejection rate of the
Metropolis-Hastings step decreases as $N$ increases.  In Fig.~\ref{fig
  ssep error ratio}(b), the error ratio for the products
$\sigma_{[xN]}\sigma_{[yN]}$ are shown for pairs $(x,y)$ at distances
ranging from ``infinitesimal'' (nearest neighbors) to $|x-y|=0.7$.
These results show that coupling control variates can effectively
improve the accuracy of calculations involving hard-to-estimate
quantities like spatial correlations.

%%%%%%%%%%%%%%%%%%%%%%%%%%%%%%%%%%%%%%%%%%%
%%%%%%%%%%%%%%%%%%%%%%%%%%%%%%%%%%%%%%%%%%%

\begin{figure}
\begin{center}
\includegraphics[bb=0in 0in 3in 2.25in]{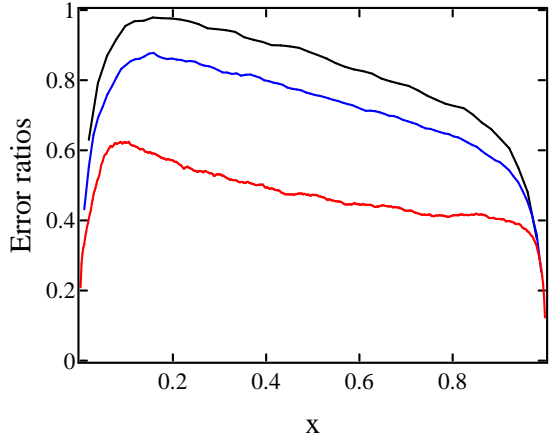}
\end{center}
\caption{The SSEP error ratio for the occupation number $\sigma_{[xN]}$
  as a function of location $x$.  The curves are, from top to bottom,
  $N=50, 100, 500$.  The parameters are $\alpha=2$, $\beta=0.1$,
  $\delta=0.3$, and $\gamma=1$.}
\label{fig ssep spatial error ratio}
\end{figure}

Fig.~\ref{fig ssep spatial error ratio} shows the error ratios for the
occupation numbers $\sigma_{[xN]}$ as as functions of spatial location
$x\in(0,1)$, for $N\in\{50,100,500\}$.  As can be seen, the error ratio
has a strong dependence on spatial location, nearly vanishing at the
boundaries but quickly attaining a near-linear profile in the interior
of the domain.  The figure show that some degrees of freedom couple
better than others, and that sites in a ``boundary layer'' near the
reservoirs couple especially well.  An explanation is that in order for
the two processes to couple at, say, site $1$, we need only that their
occupation numbers at site 1 agree, whereas for coupled moves to occur
in in the interior of the system requires that the occupation numbers of
two neighboring sites agree.  In any case, despite this dependence on
spatial location, overall the coupling control variate has improved the
accuracy of MCMC estimates by a factor of $\gtrsim 40\%$ for $N\approx
500$.

We note that the coupling control variate estimator can be implemented
with overhead of less than twice the running time of a single SSEP
simulation.  If we run two independent copies of SSEP simulations and
average the results, the standard error of the resulting estimate will
decrease by a factor of $1/\sqrt{2}\approx 0.7$, {\em i.e.}, a 30\%
gain.  We see that for single-site density estimates, the coupling
control variate offers a noticeable improvement over simply running more
copies of the simulation, and performs significantly better for two-site
estimates.

%%%%%%%%%%%%%%%%%%%%%%%%%%%%%%%%%%%%%%%%%%%%%%%%%%%%%%%%%%%%%%%%%%%
%%%%%%%%%%%%%%%%%%%%%%%%%%%%%%%%%%%%%%%%%%%%%%%%%%%%%%%%%%%%%%%%%%%

%%%%%%%%%%%%%%%%%%%%%%%%%%%%%%%%%%%%%%%%%%%%%%%%%%%%%%%%%%%%%%%%%%%
\begin{figure}
\begin{center}
\begin{tabular}{cc}
\includegraphics[bb=0in 0in 3in 2.25in]{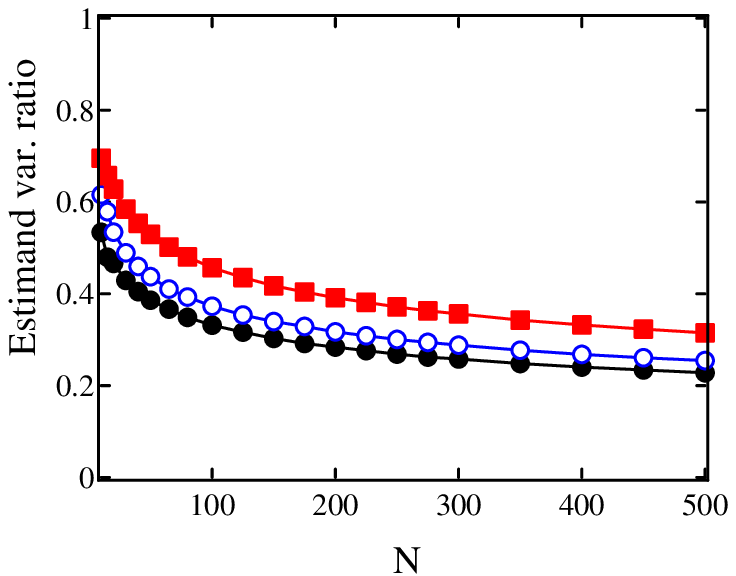}&
\includegraphics[bb=0in 0in 3in 2.25in]{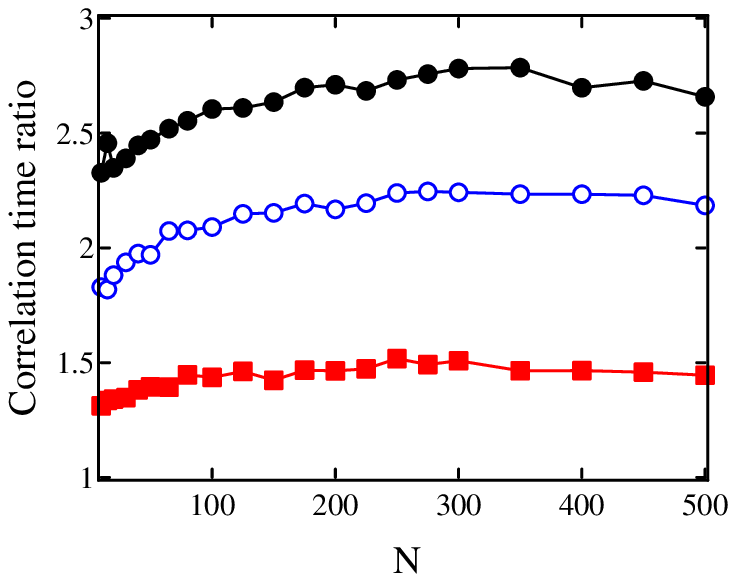}\\
(a) & (b)\\
\end{tabular}
\end{center}
\caption{The variance and correlation time components of the error ratio
  for the SSEP.  In (a), we show the factor $e_{var;N}$, as defined in
  Eq.~(\ref{error ratio components}), for $\phi=\sigma_{[xN]}$ with
  $x=0.3$ (solid discs), $x=0.5$ (open circles), and $x=0.8$ (squares).
  In (b), we show the corresponding ratios of correlation times.
  Correlation times are computed by checking numerically that Kubo
  scaling (\ref{kubo}) is in effect (batched means estimates of the
  estimator error for integration times $T\in[10^5, 10^7]$ show that the
  mean squared error $\sim T^{-1/2}$).  Then, the correlation time is
  ``reverse-engineered'' using the Kubo formula, and spot-checked by
  direct computation of time correlation functions.  Variances are
  computed by time averaging for $10^8$ time units.  The parameters are
  $\alpha=2$, $\beta=0.1$, $\delta=0.3$, and $\gamma=1$.}
\label{fig ssep var ratio}
\end{figure}
%%%%%%%%%%%%%%%%%%%%%%%%%%%%%%%%%%%%%%%%%%%%%%%%%%%%%%%%%%%%%%%%%%%

The Kubo formula (\ref{kubo}) tells us that when the simulation time $T$
is sufficiently large, the error ratio~(\ref{error ratio}) can be
written as a product of two factors:
\begin{align}
  e_N[\phi] &\approx
  \left(\frac{\var\big[\phi(\sigma)-\phi(\eta)\big]}
       {\var\big[\phi(\sigma)\big]}\right)^{1/2}
       \cdot
       \left(\frac{\tau_{couple}}{\tau}\right)^{1/2}
       \label{error ratio components}\\[1.5ex]
       &= e_{var;N}[\phi] \cdot e_{\tau;N}[\phi]~.\nonumber
\end{align}
The reasoning in Sect.~\ref{metropolization} suggests that the error
ratio $e_N$ reflects both the gain in the first factor $e_{var;N}$ by
reducing variance, and possible loss due to an increase in the second
factor $e_{\tau;N}$, by increasing correlation times.  To assess the
situation, we have plotted $e_{var,N}[\phi]$, with $\phi=\sigma_{xN}$
for a few locations $x$, in Fig.~\ref{fig ssep var ratio}(a).  This
curve should coincide with the plot of $e_N$ in Fig.~\ref{fig ssep error
  ratio}(a) {\em if} the correlation time of the SSEP were equal to that
of the coupling control variate.  Instead, we find that $e_{var;N} <
e_N$.  Fig.~\ref{fig ssep var ratio}(b) shows the ratio of integrated
autocorrelation times.  As can be seen, the coupling control variate may
increase correlation times at the same time that it reduces variance.
Here, the reduced variance wins over the increased correlation time.

%% \begin{figure}
%% \begin{center}
%% \begin{tabular}{cc}
%% \includegraphics[bb=0in 0in 3in 2.25in]{pix/ssep-prod-corrtime-ratio}\\
%% (a) & (b)\\
%% \end{tabular}
%% \end{center}
%% \caption{The SSEP correlation time ratio as a function of system size
%%   $n$.  In (a), we show the autocorrelation time ratio
%%   $\tau[\sigma_{[xn]}-\tilde\sigma_{[xn]}]/\tau[\sigma_{[xn]}]$ for the
%%   occupation number at $x=0.3$ (solid discs), $x=0.5$ (open circles),
%%   and $x=0.8$ (squares) as functions of $n$.  In (b), we show the
%%   analogous ratio for the product $\sigma_{[xn]}\sigma_{[yn]}$ with
%%   $(x,y)=(0.3,0.5)$ (discs) and for the nearest-neighbor product
%%   $\sigma_{[xn]}\sigma_{[xn]+1}$ at $x=0.6$ (squares).  The variances
%%   are computed by time averaging for $10^8$ time units.  The parameters
%%   are $\alpha=2$, $\beta=0.1$, $\delta=0.3$, and $\gamma=1$.}
%% \label{fig ssep corr time ratio}
%% \end{figure}

%{\bf [The last thing I want to add to this section: some timing
%    information on the coupling code versus the plain SSEP code, to show
%    that the coupling does not add too much computational overhead.]}

%\newpage
\subsection{KMP model}
\label{Sect KMP model}

The second model we consider is the Kipnis-Marchioro-Presutti (KMP)
model~\cite{kmp}.  This is a stochastic idealization of a chain of $N$
coupled harmonic oscillators placed at the vertices of a regular
lattice.  We think of the $i$th oscillator as having energy $\eps_i$,
given by a nonnegative real number, so that the state space is $\Omega =
[0,\infty)^N$.  Note that unlike the SSEP, $\Omega$ is uncountable.  At
  sites 0 and $N+1$, we place ``heat baths'' with temperature $T_L$ and
  $T_R$, respectively.  There are thus $N+1$ bonds in the system,
  linking site $i$ with $i\pm 1$ for $i=1, \cdots, N$.  Associated with
  each bond is an independent exponential clock of rate 1.  If the clock
  for the bond $(i,i+1)$ rings and $1\leq i\leq N-1$, then the energies
  of oscillators $i$ and $i+1$ are pooled together and redistributed
  randomly, {\em i.e.}, $\eps_i^+ = U\cdot(\eps_i^-+\eps_{i+1}^-)$ and
  $\eps_{i+1}^+ = (1-U)\cdot(\eps_i^-+\eps_{i+1}^-)$, where $U$ is a
  uniform random variable on $[0,1]$ independent of everything else,
  $\eps^+$ denotes energy after the redistribution, and $\eps^-$ denotes
  the prior energy.  If the clock for the bond $i=0$ rings, $\eps_1$
  jumps to a new energy level $u$ with probability density $\beta_L
  e^{-\beta_L u}, \beta_L = 1/T_L$.  Similarly for the bond $(N,N+1)$,
  but with parameter $\beta_R = 1/T_R$.
%% The generator for the KMP process is
%% \begin{align}
%% L_{\mbox{\tiny KMP}}f(\eps) &= \sum_{i=1}^{n-1}\int_0^1 [
%%   f(\eps^{(i,i+1),u}) - f(\eps) ]\ du\\
%% &+\int_0^\infty \beta_Le^{-\beta_Lr}[f(\eps^{1,r})-f(\eps)]\ dr
%% +\int_0^\infty \beta_Re^{-\beta_Rr}[f(\eps^{n,r})-f(\eps)]\ dr\ .\nonumber
%% \end{align}
%% where $\eps^{(i,j),u}$ is obtained from $\eps$ by setting $\eps_i =
%% u\cdot(\eps_i+\eps_j)$ and $\eps_j = (1-u)\cdot(\eps_i+\eps_j)$, and
%% $\eps^{i,r}$ from $\eps$ by setting $\eps_i=r$.
Notice that the dynamics conserves energy except at sites 1 and $N$,
just as the interior dynamics of the SSEP conserves particle number.

The KMP process provide a simple microscopic model of heat conduction.
When $T_L=T_R=T_0$, the system attains thermal equilibrium: the dynamics
satisfies detailed balance, the stationary distribution $P_N$ is a
product of Gibbs distributions with densities $\beta_0 e^{-\beta_0\eps}$
($\beta_0=1/T_0$), and the temperature at all sites is equal to $T_0$.
When $T_L\neq T_R$, we have a linear temperature profile
\begin{equation}
T(x) = T_L\cdot(1-x) + T_R\cdot x~,\qquad x\in(0,1)~,
\end{equation}
where $T(x) = \lim_{N\to\infty}\E_N[\eps_{[xn]}]$.  This non-constant
profile reflects the flow of a nonzero energy current through the
system.  The spatial correlations have a similar scaling as the
SSEP~\cite{bertini}: the limit
$$
c(x,y) =
\lim_{N\to\infty}N~\cov_N\big(\eps_{[xN]},\eps_{[yN]}\big)
$$
exists, and
$$
c(x,y) \propto (T_R-T_L)^2\cdot x(1-y)\ ,\qquad 0 < x < y < 1\ .
$$
Like the SSEP, $\cov_N(\eps_{[xN]},\eps_{[yN]}) = O(1/N)$.  Thus, one
encounters similar difficulties when estimating spatial correlations
numerically.

\begin{figure}
\begin{center}
\input{pix/kmp-couple.tex}
\end{center}
\caption{Illustration of the KMP coupling.  Because the interaction
  conserves energy, the point $(X_i,X_{i+1})$ is constrained to lie on
  the line $X_i + X_{i+1} =\mbox{const}$ both before and after the
  interaction.}
\label{fig kmp couple}
\end{figure}
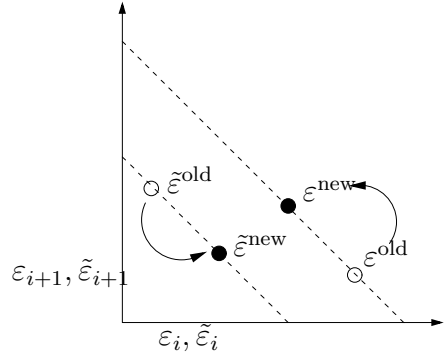

\bigskip
It has been shown that the KMP model attains LTE as $N\to\infty$, {\em
  i.e.}  $k$-site marginals converge to a product of Gibbs
distributions, with a local temperature $T(x)$ given by the linear
profile above.  This suggests that we use
\begin{equation}
Q_N(\eps) = \prod_{i=1}^N \beta_ie^{-\beta_i\eps_i}~,
\end{equation}
where $\beta_i = 1/T(x_i)$, as approximate stationary distribution.  A
simple coupling of the KMP process to itself is also available: given
two copies of the KMP process, we make the same bonds ``ring'' at the
same time.  For interior bonds, we use the same uniform random numbers
$U$ to split energy in both copies; for heat baths, we set the boundary
sites to the same new energy.
%% : its generator is
%% \begin{align}
%% L_{\mbox{\tiny KMP-couple}}f(\eps,\tilde\eps) &= \sum_{i=1}^{n-1}\int_0^1 du\ [
%%   f(\eps^{(i,i+1),u},\tilde\eps^{(i,i+1),u}) - f(\eps,\tilde\eps) ]\\
%% &+\int_0^\infty dr\ \beta_Le^{-\beta_Lr}[f(\eps^{1,r},\tilde\eps^{1,r})-f(\eps,\tilde\eps)]\nonumber\\
%% &+\int_0^\infty dr\ \beta_Re^{-\beta_Rr}[f(\eps^{n,r},\tilde\eps^{n,r})-f(\eps,\tilde\eps)]\ .\nonumber
%% \end{align}
The coupling is illustrated in Fig.~\ref{fig kmp couple}: it entails
having the $\tilde\eps$ process use the same ``randomness'' as the
$\eps$ process to redistribute energy between nearby sites.

%%%%%%%%%%%%%%%%%%%%%%%%%%%%%%%%%%%%%%%%%%%%%%%%%%%%%%%%%%%%%%%%%%%
%%%%%%%%%%%%%%%%%%%%%%%%%%%%%%%%%%%%%%%%%%%%%%%%%%%%%%%%%%%%%%%%%%%

One difference from the SSEP is that the KMP model has an uncountable
state space, so the algorithm described in Sect.~\ref{sect1: coupling
  control variates} requires slight modification.  This is
straightforward for Markov jump processes with transition densities: one
can simply replace the ratio of transition rate coefficients $R_X$ in
Eq.~(\ref{eq metropolis ratio}) with the ratio of the corresponding
densities.  The KMP process does not only have an uncountable state
space, though --- it also has singular transition rate measures (this is
a consequence of energy conservation).  Nonetheless, it can be checked
that the ratios are well-defined in this case, and yield the following
Metropolis ratios:
\begin{center}
\begin{tabular}{rl}
Interaction resulting in $(\eps_i,\eps_j)\mapsto(\eps'_i,\eps'_j)$, $|i-j|=1$
& $Z_{ij} = \exp\Big([\beta_i\eps_i +
  \beta_j\eps_j] - [\beta_i\eps'_i + \beta_j\eps'_j]\Big)$\\
Left heat bath setting $\eps_1\mapsto\eps'_1$
& $Z_L=\exp\big((\beta_L-\beta_1)\cdot(\eps'_1-\eps_1)\big)$\\
Right heat bath setting $\eps_n\mapsto\eps'_n$
& $Z_R=\exp\big((\beta_R-\beta_n)\cdot(\eps'_n-\eps_n)\big)$\\
\end{tabular}
\end{center}
Applying the algorithm in Sect.~\ref{metropolization} with these ratios
yields a coupling control variate which preserves the local equilibrium
distribution $Q_N$.

%%%%%%%%%%%%%%%%%%%%%%%%%%%%%%%%%%%%%%%%%%%%%%%%%%%%%%%%%%%%%%%%%%%
\begin{figure}
\begin{center}
\begin{tabular}{cc}
\includegraphics[bb=0in 0in 3in 2.25in]{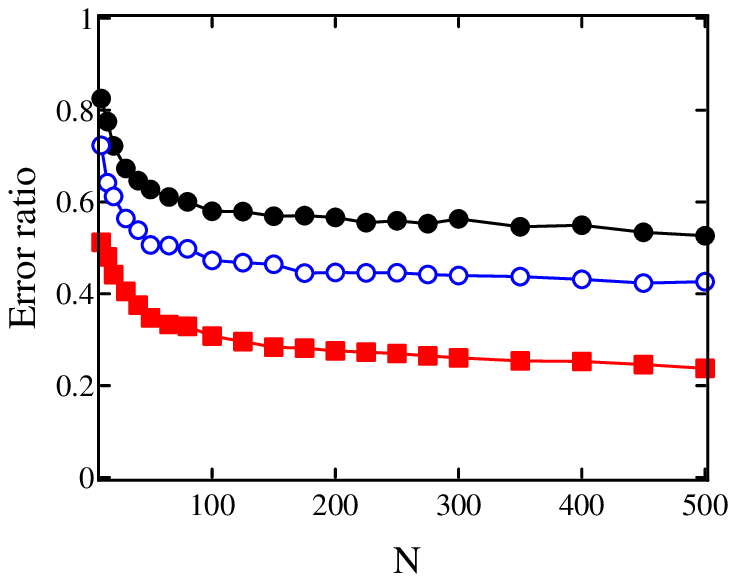}&
\includegraphics[bb=0in 0in 3in 2.25in]{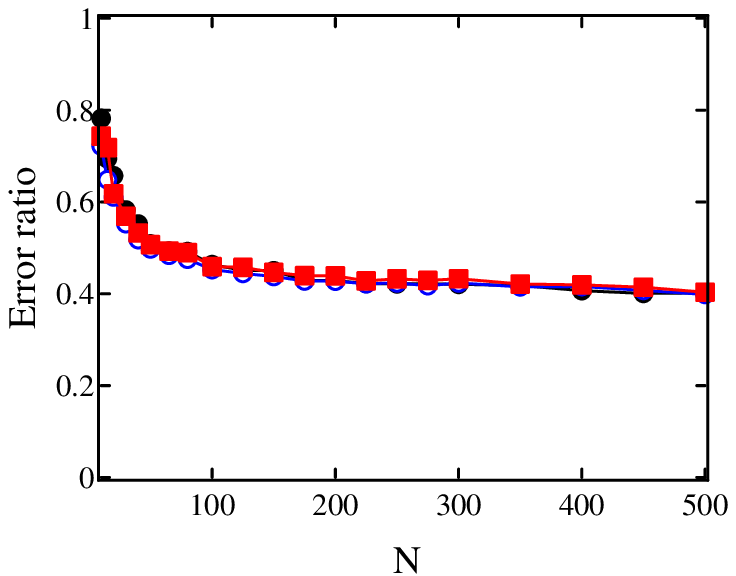}\\
(a) & (b)\\
\end{tabular}
\end{center}
\caption{The KMP error ratio as a function of system size $N$.  In (a),
  we show the error ratio for the estimated mean energies at $x=0.3$
  (solid discs), $x=0.5$ (open circles), and $x=0.8$ (squares) as
  functions of $N$.  In (b), we show the error ratio for the
  near-neighbor product $\eps_{[xN]}\cdot\eps_{[xN]+1}$ with $x=0.5$
  (solid discs), and for the products $\eps_{[xN]}\cdot\eps_{[yN]}$ with
  $(x,y)=(0.4,0.7)$ (open circles) and $(x,y)=(0.2,0.9)$ (squares).  The
  errors are estimated using batched means estimator.
  %% The errors are estimated by subdividing a simulation of duration
  %% $10^9$ time units into $10^4$ batches of duration $10^5$ each and
  %% forming a batched means estimate.
  The parameters are $T_L=10$ and $T_R=100$.}
\label{fig kmp error ratio}
\end{figure}
%%%%%%%%%%%%%%%%%%%%%%%%%%%%%%%%%%%%%%%%%%%%%%%%%%%%%%%%%%%%%%%%%%%

%%%%%%%%%%%%%%%%%%%%%%%%%%%%%%%%%%%%%%%%%%%%%%%%%%%%%%%%%%%%%%%%%%%
\begin{figure}
\begin{center}
\begin{tabular}{cc}
\includegraphics[bb=0in 0in 3in 2.25in]{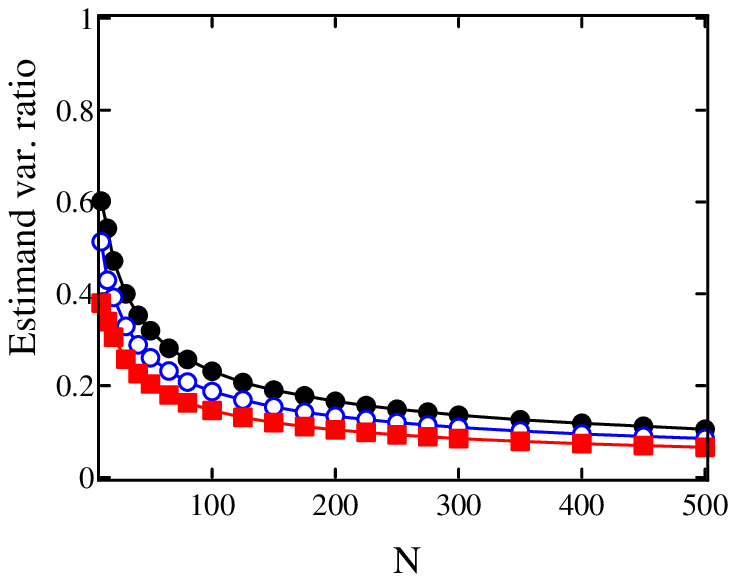}&
\includegraphics[bb=0in 0in 3in 2.25in]{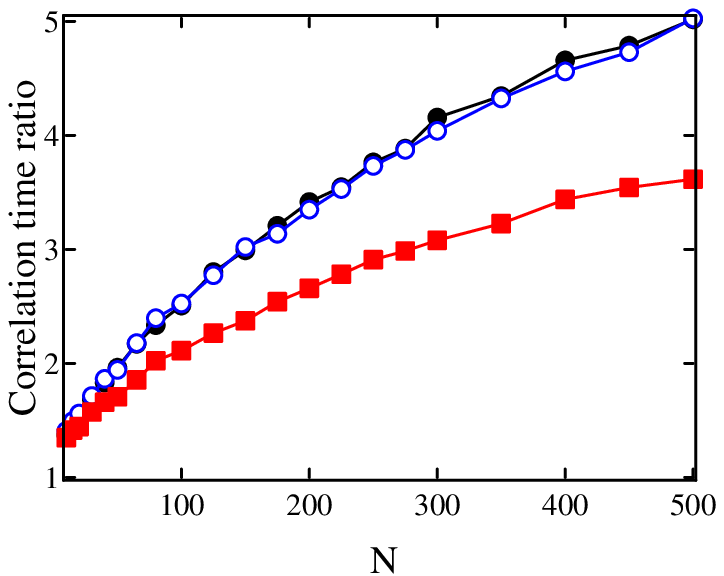}\\
(a) & (b)\\
\end{tabular}
\end{center}
\caption{The variance and correlation time components of the error ratio
  for the KMP process.  In (a), we show the factor $e_{var;N}$, as
  defined in Eq.~(\ref{error ratio components}), for $\phi=\eps_{[xN]}$
  for $x=0.3$ (solid discs), $x=0.5$ (open circles), and $x=0.8$
  (squares).  In (b), we show the corresponding ``reverse-engineered''
  ratio of correlation times.  The parameters are $T_L=10$ and
  $T_R=100$.}
\label{fig kmp var ratio}
\end{figure}
%%%%%%%%%%%%%%%%%%%%%%%%%%%%%%%%%%%%%%%%%%%%%%%%%%%%%%%%%%%%%%%%%%%

\bigskip
\noindent
{\em Numerical results.}  Fig.~\ref{fig kmp error ratio}(a) shows the
error ratios for various sites in the KMP model.  As is the case for the
SSEP, the coupling control variate significantly reduces the variance of
the estimator.  In contrast to the SSEP, the amount by which the error
is reduced depends more strongly on location, ranging from 20 to 60\%.
Fig.~\ref{fig kmp var ratio}(b) shows the error ratios for the products
$\eps_{[xN]}\cdot\eps_{[yN]}$ for pairs $(x,y)$ located at various
distances.  These ratios are much more consistent and tend to $\approx
40\%$ for the range of $N$ tested.

Fig.~\ref{fig kmp var ratio}(a) shows the corresponding factor
$e_{var;N}$.  As in the case of the SSEP, $e_{var;N}$ is strictly
smaller than the error ratio $e_N$; at the same time, the ratio
$e_{\tau;N}$ of correlation times increase; see Fig.~\ref{fig kmp var
  ratio}(b).  Thus, Metropolis rejections can have a dramatic effect on
the correlation time of the coupling control variate.  Despite that, the
overall performance of the coupling control variate estimator is quite
good: even at its worst, the accuracy has been improved by 40\%.

%%%%%%%%%%%%%%%%%%%%%%%%%%%%%%%%%%%%%%%%%%%
%\newpage
\section*{Conclusion}

We have shown that Markov couplings, when available, can be used
effectively to improve the accuracy of Markov chain Monte Carlo
calculations.  This method useful in situations where the stationary
distribution is not known explicitly, as in the case of nonequilibrium
transport models.  As shown by the examples considered in this paper,
good candidates for approximate stationary distribution can be found
based on physical reasoning, and when an effective coupling is available
for the Markov process at hand, one can construct an effective coupling
control variate.

The numerical results suggest various directions for improvement.  In
particular, the observation that coupling control variate has larger
correlation times than the original process suggests that one try to
``trade'' variance for correlation time.  However, simple ideas like
resampling the energy of random sites at random times, as in heat bath /
partial resampling, may very well increase variance more than it
decreases correlation time, resulting in a net gain of error.  A related
issue is the dependence of the estimator error ratio on observables: in
many applications, it is desirable to be able to optimize the error
ratio only for observables of interest.  (One does not expect to be able
to have small error ratios for all observables unless the approximate and
true stationary distributions are close in the total variation norm.)

Finally, we mention that it might be possible to use related coupling
methods for sensitivity analysis.  If the Markov process depends on
parameters $\theta$, then the observable $\phi$ in Eq.~(\ref{ccv})
becomes $\phi_\theta$ and the sensitivities are derivatives of
$\phi_\theta$ with respect to $\theta$.  Sensitivities are used, for
example, in numerical computation of optimal stochastic controls in
situations where the curse of dimensionality makes dynamic programming
impractical.  If there is a known formula for the stationary
distribution $P_\theta$, two common methods for evaluating
sensitivities, the {\em common random variables} (or {\em same paths})
method\footnote{An infinitesimal variation version of this method
  sometimes is associated with the Malliavin calculus.}  and the {\em
  likelihood ratio} (or {\em score function}) methods.
Glynn~\cite{glynn} and others have generalizations of the likelihood
ratio method to situations where $T$ is known but not $P$.  It also
might be helpful to have such a generalization of the same paths method.

\paragraph{Acknowledgments.}
We thank Peter Glynn for pointing out the references~\cite{glynn,o1,o2}.
JG was supported in part by the Applied Mathematical Sciences Program of
the U.S. Department of Energy under Contract DEFG0200ER25053.  KL was
supported in part by an NSF Mathematical Sciences Postdoctoral
Fellowship.

%%%%%%%%%%%%%%%%%%%%%%%%%%%%%%%%%%%%%%%%%%%
%\newpage

\end{document}

%% file: pix/kmp-couple.tex
\setlength{\unitlength}{0.0005in}
\begingroup\makeatletter\ifx\SetFigFont\undefined%
\gdef\SetFigFont#1#2#3#4#5{%
  \reset@font\fontsize{#1}{#2pt}%
  \fontfamily{#3}\fontseries{#4}\fontshape{#5}%
  \selectfont}%
\fi\endgroup%
{\renewcommand{\dashlinestretch}{30}
\begin{picture}(3762,3627)(0,-10)
\put(1605,990){\makebox(0,0)[lb]{{\SetFigFont{12}{14.4}{\familydefault}{\mddefault}{\updefault}$\tilde\eps^{\rm new}$}}}
\put(1065.000,1368.214){\arc{829.164}{0.9971}{3.5821}}
\blacken\thicklines
\path(1192.745,943.570)(1290.000,1020.000)(1168.220,998.329)(1192.745,943.570)
\thinlines
\put(2175,1515){\blacken\ellipse{150}{150}}
%\put(2175,1515){\ellipse{150}{150}}
\put(1455,1020){\blacken\ellipse{150}{150}}
%\put(1455,1020){\ellipse{150}{150}}
\put(2865,795){\whiten\ellipse{150}{150}}
%\put(2865,795){\ellipse{150}{150}}
\put(750,1695){\whiten\ellipse{150}{150}}
%\put(750,1695){\ellipse{150}{150}}
\path(450,300)(450,3600)
\blacken\thicklines
\path(472.500,3525.000)(450.000,3600.000)(427.500,3525.000)(472.500,3525.000)
\thinlines
\path(450,300)(3750,300)
\blacken\thicklines
\path(3675.000,277.500)(3750.000,300.000)(3675.000,322.500)(3675.000,277.500)
\thinlines
\dashline{60.000}(450,3225)(3375,300)
\drawline(3375,300)(3375,300)
\dashline{60.000}(450,2025)(2175,300)
\put(825,0){\makebox(0,0)[lb]{{\SetFigFont{12}{14.4}{\rmdefault}{\mddefault}{\updefault}$\eps_i,
\tilde\eps_i$}}}
\put(-700,675){\makebox(0,0)[lb]{{\SetFigFont{12}{14.4}{\rmdefault}{\mddefault}{\updefault}$\eps_{i+1}, \tilde\eps_{i+1}$}}}
\put(2955,885){\makebox(0,0)[lb]{{\SetFigFont{12}{14.4}{\familydefault}{\mddefault}{\updefault}$\eps^{\rm old}$}}}
\put(2325,1560){\makebox(0,0)[lb]{{\SetFigFont{12}{14.4}{\familydefault}{\mddefault}{\updefault}$\eps^{\rm new}$}}}
\put(915,1650){\makebox(0,0)[lb]{{\SetFigFont{12}{14.4}{\familydefault}{\mddefault}{\updefault}$\tilde\eps^{\rm old}$}}}
\put(2917.500,1350.000){\arc{762.053}{4.5343}{6.9149}}
\blacken\thicklines
\path(2969.632,1756.436)(2850.000,1725.000)(2970.351,1696.440)(2969.632,1756.436)
\end{picture}
}